\documentclass[journal]{new-aiaa}
\usepackage[utf8]{inputenc}
\usepackage{textcomp}
\usepackage[super]{nth}
\usepackage{graphicx}
\usepackage{amsmath}
\usepackage[version=4]{mhchem}
\usepackage{siunitx}
\usepackage{longtable,tabularx}
\usepackage[normalem]{ulem}
\newcommand\redout{\bgroup\markoverwith
{\textcolor{red}{\rule[.5ex]{2pt}{0.4pt}}}\ULon}
\usepackage{xcolor}

\usepackage{hyperref} 
\hypersetup{
    colorlinks = true,  
	allcolors = black  
}
\usepackage{multicol}
\usepackage{xtab} 
\usepackage{bm}
\RequirePackage{color}
\usepackage{xcolor}
\usepackage{soul}

\usepackage{comment}
\usepackage{subfigure}

\renewcommand{\arraystretch}{0.6}

\def\bfr{\mbox{\boldmath $r$}}

\def\bfv{\mbox{\boldmath $v$}}

\def\bfx{\mbox{\boldmath $x$}}

\newcommand{\norm}[1]{\left\lVert#1\right\rVert}

\title{Convex Optimization of Linear Impulsive Rendezvous}

\author{Boris Benedikter\footnote{PhD Student, Department of Mechanical and 
Aerospace Engineering, Sapienza University of Rome, Via Eudossiana 18 - 00184, 
Rome, Italy; boris.benedikter@uniroma1.it}
and 
Alessandro Zavoli\footnote{Research Assistant, Department of Mechanical and 
Aerospace Engineering, Sapienza University of Rome, Via Eudossiana 18 - 00184, 
Rome, Italy; alessandro.zavoli@uniroma1.it}}
\affil{Sapienza University of Rome, Rome, 00184, Italy}

\begin{document}

\maketitle

\section*{Nomenclature}

\medskip

{\renewcommand\arraystretch{1.0}
\noindent\begin{tabular}{l l p{0.39\linewidth} l l p{0.39\linewidth}}
$a$ & $=$ & semi-major axis & $\bm{\Delta v}$ & $=$ & velocity change vector \\
$B$ & $=$ & actuation matrix & $\theta$ & $=$ & true anomaly \\
$e$ & $=$ & eccentricity & $\rho$ & $=$ & $ 1 + e \cos\theta$ \\
$J$ & $=$ & objective function & $\sigma$ & $=$ & free auxiliary variable \\
$k^2$ & $=$ & $ h / p^2$ & $\Phi$ & $=$ & state transition matrix \\
$I$ & $=$ & identity matrix & $\Omega$ & $=$ & right ascension of ascending node \\
$i$ & $=$ & inclination & $\omega$ & $=$ & argument of periapsis \\
$h$ & $=$ & angular momentum & $\norm{\cdot}$ & $=$ & $L^2$ norm \\
$M$ & $=$ & number of grid nodes & \multicolumn{2}{@{}l}{\textbf{Subscripts}}\\
$N$ & $=$ & number of impulses & $0$ & $=$ & initial boundary \\
$\bm{r}$ & $=$ & chaser relative position & $f$ & $=$ & final boundary \\
& & in the target orbital frame $(x, y, z)$ & \multicolumn{2}{@{}l}{\textbf{Superscripts}}\\
$p$ & $=$ & semilatus rectum & $\sim$ & $=$ & transformed variable \\
$t$ & $=$ & time & $\prime$ & $=$ & derivative with respect to $\theta$ \\
$\bm{v}$ & $=$ & chaser relative velocity & $-$ & $=$ & values immediately before an impulse \\
& & in the target orbital frame $(v_x, v_y, v_z)$ & $+$ & $=$ & values immediately after an impulse \\
$\bm{x}$ & $=$ & state vector \\
\end{tabular}}

\section{Introduction}

\lettrine{T}{his} note aims at providing a concise and self-contained document that 
describes a clear and easy-to-understand method, that could be useful for 
a reader that is approaching the linear-impulsive rendezvous topic for the 
first time, but that is also flexible enough to accommodate problem variations,
such as additional constraints like bounded-magnitude $\Delta v$, or direction 
constraints with minor modifications.

A convex approach for the optimization of 
time-fixed non-cooperative rendezvous problems in a linear relative dynamic field is proposed. 
Despite its simplicity, linear dynamics is very effective
to describe the relative motion between two spacecraft (or between one 
spacecraft and a reference, virtual, satellite), and it is routinely employed
in several practical scenarios, such as 
spacecraft docking~\cite{weiss2015model}, 
proximity operations~\cite{subbarao2008nonlinear}, 
formation flying~\cite{vassar1985formationkeeping}, 
collision avoidance~\cite{epenoy2011fuel},
and even in unconventional scenarios~\cite{zavoli2019gtoc}.
Starting with the Clohessy–Wiltshire or Hill 
equations~\cite{hill1878researches},
that are only valid for circular reference orbits,
numerous improved versions have been developed in order to account for
large-angle gaps~\cite{Baranov1990},
elliptic orbit~\cite{tschauner1965rendezvous, Yamanaka2002},
second-order terms~\cite{Kechichian1992}, 
and even orbital perturbations~\cite{Schweighart2002}.
In this note, the Tschauner-Hempel 
equations~\cite{tschauner1965rendezvous},  
valid for reference elliptic orbits of arbitrary eccentricity,
are adopted.

The rendezvous problem, 
that is, the problem of an active, \emph{chaser},  spacecraft that must reach 
for a passive, \emph{target}, spacecraft in a given amount of time is a 
well-known topic in spaceflight mechanics.
A number of solution approaches have been developed for both far-field 
rendezvous, involving a nonlinear two-body dynamics in either
passive~\cite{Mirfakhraie1994,lu2013autonomous} or cooperative 
scenarios~\cite{benedikter2019convexrendezvous, Zavoli2015}, 
as well as for the case of the linear dynamics case here investigated.
In the latter case, most of the published works 
focused on the use of the Pontryagin Maximum Principle (PMP) for deriving  
optimality conditions for propellant-optimal trajectories. %
The pioneering work of Neustadt~\cite{neustadt1964optimization} states that 
for linear dynamics in an $n$-dimensional state-space
the optimal transfer requires at most $n$ impulses.
By applying the primer vector theory, Lawden derived a set of first-order 
necessary conditions~\cite{lawden1963general}. 
Subsequent efforts due to Prussing~\cite{prussing1994optimal}, 
Carter~\cite{carter1991optimal}, and Jezewski~\cite{Jezewski1980} 
led to the definition of sufficient conditions for the optimality
of a coplanar linear rendezvous problem.
Leveraging on these conditions, analytical methods were developed for 
calculating burn time, magnitude, and direction of the impulsive maneuvers for 
a fixed number of 
burns~\cite{prussing1969optimal,carter2000quadratic}.
However, the optimal number of impulses cannot be predicted in advance, and 
it is usually adjusted iteratively, by inspecting the solution and adding 
maneuvers or coasting arcs as needed according to PMP~\cite{Lion1968}.
More recently, general closed-form solutions have been derived for linear 
impulsive rendezvous using accurate dynamics based on relative orbit 
elements~\cite{chernick2017new}.


In this work, 
the original optimal control problem is transformed into a 
convex one  by discretizing the independent variable  over a  (sufficiently dense) grid and constraining the impulses to be located at the grid nodes.
As a result of the \emph{convexification}, one has a  
theoretical guarantee of convergence towards the global optimum in a 
limited, short, time regardless of the initialization.
The optimal number of impulses and their epochs (up to some discretization
precision) are directly obtained, 
without the need for any \emph{a priori} assumptions on the solution structure. 
This is quite a desirable feature of the proposed approach, as the problem 
finds many applications in time-critical scenarios where both the 
computational efficiency and the reliability of the algorithm are 
primary requirements.

This manuscript is organized as follows.
Section~\ref{sec:problem_description} describes the original optimal control
problem.
In Section~\ref{sec:convex_formulation}, the problem is transcribed into a 
convex optimization problem.
Numerical results are presented in Section~\ref{sec:numerical_results} and compared with those provided by other solution methods.
A conclusion section ends the note.

%
\section{Problem Description} \label{sec:problem_description}
\subsection{System Dynamics}

\begin{figure}[hbt!]
    \centering
    \includegraphics[width=.5\textwidth]{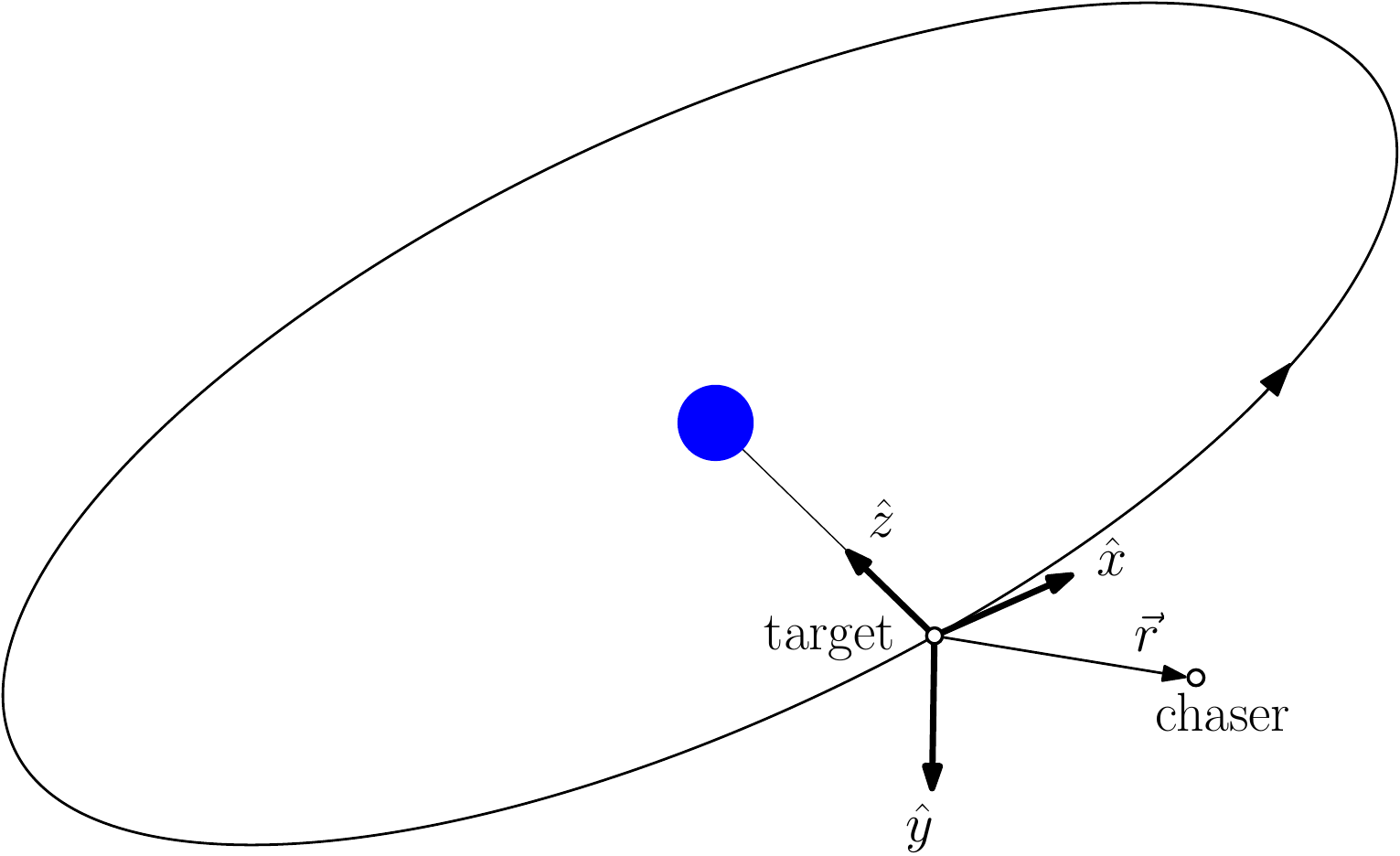}
    \caption{Orbital Reference frame 
    }
    \label{fig:reference_frame}
\end{figure}

This section introduces the Tschauner-Hempel equations that allow describing 
the linear relative dynamics of a spacecraft
with respect to a reference point that moves on a keplerian elliptic orbit of 
arbitrary eccentricity. 
The chaser state $\bfx$ is completely described by its relative position
$\bfr$ and relative velocity $\bfv$ with respect to the target, that is,
$\bm{x} = \begin{bmatrix}
    \bm{r} & \bm{v}
\end{bmatrix}^T$,
that are conveniently expressed in a rotating reference frame, 
centered in the target spacecraft with the $z$-axis in 
the radial direction and pointing towards the Earth and the
$y$-axis in the direction opposite to the angular velocity vector of the target,
as shown in Figure~\ref{fig:reference_frame}.

Let $\tilde\bfx = \begin{bmatrix}\tilde\bfr & \tilde\bfv \end{bmatrix}$ 
be the transformed state, obtained by applying the ``direct'' transformation:
\begin{align}
    \bm{\tilde{r}} = \rho \bm{r} && 
    \bm{\tilde{v}} = -e \sin\theta \bm{r} + 1 / (k^2 \rho)\bm{v}
\end{align}
where $e$ is the eccentricity of the target's orbit, 
and $\rho$ and $k$ are defined as
\begin{gather}
    \rho = 1 + e \cos\theta \\
    k^2 = h / p^2
\end{gather}
where $h$ and $p$ are the target's orbit angular momentum and semilatus rectum,
respectively.
The corresponding inverse transformation is readily available:
\begin{align}
\label{eq:inverseTransformation}
    \bm{r} = 1 / \rho \bm{\tilde{r}} &&
    \bm{v} = k^2 (e \sin\theta \bm{\tilde{r}} + \rho \bm{\tilde{v}})
\end{align}

Thanks to the coordinate transformation adopted, 
the equations of relative motion become:
\begin{align}
    \tilde{x}^{\prime\prime} &= 2 \tilde{z}^\prime \\
    \tilde{y}^{\prime\prime} &= -\tilde{y} \\
    \tilde{z}^{\prime\prime} &= 3 \tilde{z} / \rho -2\tilde{x}^\prime    
\end{align}
where $( \ )^\prime$ denotes the derivative with respect to the true anomaly 
$\theta$, which is used as independent variable instead of time $t$.
Please notice that this is a non-autonomous system as $\rho$ 
depends on the independent variable $\theta$.

Alternatively, the chaser dynamics can be effectively described in terms of 
the state transition matrix (STM) $\Phi(t_0, t)$ 
that links the spacecraft transformed state $\bm{\tilde{x}}(t)$ at time $t$, 
with the transformed state $\bm{\tilde{x}}(t_0)$ at time $t_0$:
\begin{equation}
    \bm{\tilde{x}}(\theta) = \Phi(\theta, \theta_0) \bm{\tilde{x}}(\theta_0)
    \label{eq:state_transition_matrix_anu}
\end{equation}

The STM $\Phi(\theta, \theta_0)$ can be decomposed into
STMs for in-plane and out-of-plane motions, respectively.
The in-plane relative state at arbitrary $\theta$ can be calculated as: 
\begin{equation}
    \begin{bmatrix}
        \tilde{x} \\
        \tilde{z} \\
        \tilde{v}_{x} \\
        \tilde{v}_{z}                
    \end{bmatrix} = 
    \Phi_\theta \Phi^{-1}_{\theta_0}
    \begin{bmatrix}
        \tilde{x}_0 \\
        \tilde{z}_0 \\
        \tilde{v}_{x_0} \\
        \tilde{v}_{z_0}        
    \end{bmatrix}
\end{equation}
where:
\begin{align}
    \Phi_\theta &=
    \begin{bmatrix}
        1 & -c(1 + 1 / \rho) & s(1 + 1 / \rho) & 3 \rho^2 J \\
        0 & s & c & 2 - 3 esJ \\
        0 & 2s & 2c - e & 3 (1 - 2esJ) \\
        0 & s' & c' & -3e(s' J + s / \rho^2)
    \end{bmatrix}_\theta \\    
    \Phi^{-1}_{\theta_0} &= \frac{1}{1 - e^2} 
    \begin{bmatrix}
        1 - e^2 & 3es (1 / \rho + 1 / \rho^2) & - es(1 + 1 / \rho) & -ec + 2 \\
        0 & -3s(1 / \rho + e^2 / \rho^2) & s(1 + 1 / \rho) & c - 2e \\
        0 & -3(c / \rho + e) & c(1 + 1 / \rho) & -s \\
        0 & 3 \rho + e^2 - 1 & -\rho^2 & es
    \end{bmatrix}_{\theta_0}
\end{align}
and
\begin{align}
    s &= \rho \sin\theta \\
    c &= \rho \cos\theta \\
    s' &= \cos\theta + e \cos(2 \theta) \\
    c' &= -(\sin\theta + e \sin(2 \theta)) \\
    J &= k^2 (t - t_0)
\end{align}
$s'$ and $c'$ denote the derivatives of $s$ and $c$ with respect to the true
anomaly $\theta$.
Instead, $J$ is computed using time $t$, 
that is related to the true anomaly by Kepler's equation.
Finally, the out-of-plane components of the relative state 
can be computed as:
\begin{equation}
    \begin{bmatrix}
        \tilde{y} \\
        \tilde{v}_y
    \end{bmatrix} = 
    \frac{1}{\left.\rho\right|_{{\theta - \theta_0}}}
    \begin{bmatrix}
        c & s \\
        -s & c 
    \end{bmatrix}_{\theta - \theta_0}
    \begin{bmatrix}
        \tilde{y}_0 \\
        \tilde{v}_{y_0}
    \end{bmatrix}
\end{equation}

\subsection{Optimal Control Problem}

An impulsive thrust model is assumed, that is,
instantaneous velocity changes 
$\bm{\Delta v_j} = 
\left[ \Delta v_{x,j} \; \Delta v_{y,j} \; \Delta v_{z,j} \right]$ 
are applied at some, unknown, anomalies $\theta_j$, for $j = 1, \ldots, N$. 
By introducing the compact notation $\bm{x}_{j} = \bm{x}(\theta_{j})$, 
and using 
superscripts ``-'' and ``+'' to identify the state values immediately before 
and after the burns, respectively, the general optimization problem can be 
defined in the compact form:
\begin{align}
    \min_{N,\, \bm{\Delta \tilde{v}}_j,\, \theta_j}&
    J = \sum_{i=1}^{N}{\Delta v_j} &
    \label{eq:objective_original} \\
    \text{subject to} 
        &\bm{\tilde{x}}_1^- = \bm{\tilde{x}}_0 & 
        \label{eq:initial_conditions_original}\\
        &\bm{\tilde{x}}_{N}^+ = \bm{\tilde{x}}_f & 
        \label{eq:final_conditions_original}\\
        &\bm{\tilde{x}}_{j+1}^- = \Phi(\theta_{j+1},\theta_j) \, 
        \bm{\tilde{x}}_j^+ 
        &&j = 1, \dots, N-1 
        \label{eq:differential_constraints_original} \\
        &\bm{\tilde{x}}_{j}^+ = \bm{\tilde{x}}_j^- 
        + B \, \bm{\Delta \tilde{v}}_j 
        &&j = 1, \dots, N
        \label{eq:before_after_burns_original}
\end{align}
where $\Delta v = \norm{\bm{\Delta v}}$ is the magnitude of the velocity change 
and $B = [0_{3\times 3} \, I_{3\times 3}]^T$ is the actuation matrix.
Equations~\eqref{eq:initial_conditions_original} 
and \eqref{eq:final_conditions_original} 
are the initial and terminal conditions, respectively.
The initial and final states, $\bm{\tilde{x}}_0$ and $\bm{\tilde{x}}_f$,
are assigned and depend on the specific problem instance. 
As an example, by enforcing $\bm{\tilde{x}}_f = \bm{0}$ 
the rendezvous condition is attained. 
Nevertheless, this formulation allows for arbitrary terminal conditions 
and it is also suitable for formation reconfiguration.

Please notice that not only  the locations of the velocity impulses are 
unknown,
but also the optimal number of burns $N$, even though it is
constrained in $\mathcal{I} = [1, N_{max}]$,
where $N_{max}$ is the cardinality of the state space, that is, 
4 in a planar problem, or 6 in a three-dimensional case, due to Neustadt 
necessary optimality conditions~\cite{neustadt1964optimization}.  

\section{Convex Formulation}
\label{sec:convex_formulation}

The optimal control problem formulated in 
Equations~\eqref{eq:objective_original}--\eqref{eq:before_after_burns_original}
is a general, nonlinear optimization problem since the STM is a nonlinear 
function of the unknown impulse locations $\theta_j$.
Moreover, the problem depends on the optimal number of impulses, 
which is not known in advance.
In general, such a problem might be hard to solve.
A common approach consists of fixing $N$ and 
searching for the location of the impulse, assuming that
the first and the last $\Delta v$ to be applied 
at the initial and final time, respectively.
In particular, $N$ can be set equal to the upper bound $N_{max}$ on the number of impulses. 
However, optimal solutions of typical problems may require fewer impulses~\cite{prussing2003optimal}, therefore such an approach may not be suitable.

In this work, we propose a numerical method that requires no \emph{a priori} decision 
on the number of impulses nor on their location.
Our approach converts the presented optimal control problem
into a Second-Order Cone Programming (SOCP) problem, a special
class of convex optimization problems, characterized by 
a linear objective function, linear equality constraints,
and second-order cone constraints. 
This class of programming problems allows for representing quite complex 
constraints and can be solved by means of highly-efficient interior point 
methods, even for a large number of variables~\cite{alizadeh2003second}.
The convexification is carried out by discretizing the problem over a 
sufficiently dense $\theta$-grid and by constraining the impulses to 
be located at the grid nodes.
In this way, the locations and the total number of impulses are no longer 
optimization variables, and the corresponding nonlinearities are avoided.

\subsection{Discretization Grid}
The independent variable, here the target's true anomaly $\theta$,
is discretized into a grid of $M$ points, also referred to as mesh nodes:
\begin{equation}
    \theta_0 = \theta_1 < \dots < \theta_{M} = \theta_f
    \label{eq:nu_mesh}    
\end{equation}

The finite discretization of the independent variable introduces 
a discretization error on the locations of the impulses, 
as an instantaneous velocity change can be applied only at those points.
State and control variables are discretized over the chosen grid, 
the differential constraints for each interval can be transformed into 
algebraic constraints involving the STM:
\begin{equation}
    \bm{\tilde{x}}_{j + 1}^{-} = 
    \Phi(\theta_{j + 1}, \theta_{j}) \, \bm{\tilde{x}}_j^{+} 
    \qquad \qquad j = 1, \dots, M - 1 \\
    \label{eq:differential_con_discrete}
\end{equation}

Even though the proposed approach poses no restriction on the grid spacing, 
a uniform node distribution is considered hereafter, as usually one has no 
prior information on the optimal location of the maneuvers. 
Intuitively, the quality of the approximation is related to the distance
between the nodes and the (true) optimal locations of the impulses:  
a denser grid will probably allow for a more accurate solution, 
but will also lead to a more computationally expensive problem. 
A parametric analysis will be conducted to investigate the effects of 
the number of discretization points on the overall quality of the attained
numerical results.

\subsection{Objective Function}
The goal of the optimization is to minimize the overall $\Delta v$.
Thus, the objective function $J$ can be expressed as the sum of the 
$\text{L}^2$ norm of the velocity change vectors:
\begin{equation}
    J = \sum_{j = 1}^{M} \norm{
    \bm{\Delta {v_j}}}_2
    = \sum_{j = 1}^{M}{k^2 \rho_j \norm{\bm{\Delta \tilde{v}}_j}_2}
    \label{eq:objective_norm2}
\end{equation}
where Eq.~\eqref{eq:inverseTransformation}
has been used to relate the transformed and actual velocity change, and 
$\rho_j = 1 + e\cos(\theta_j)$ are $M$ known constants, evaluated at the grid nodes. 

This formulation provides a convex objective function.
However, in a SOCP problem, the objective function must be a linear function
of the optimization variables.
Hence, we introduce the auxiliary variables $\sigma$ 
that represent a free upper bound to the 
$\text{L}^2$ norm of the velocity change vectors as:
\begin{equation}
    \Delta \tilde{v}_{x, j}^2 + \Delta \tilde{v}_{y, j}^2 
    + \Delta \tilde{v}_{z, j}^2 \leq \sigma_j^2
    \qquad \qquad j = 1, \dots, M
    \label{eq:cone_constraint}
\end{equation}
Notice that Eq.~\eqref{eq:cone_constraint} is a second-order cone constraint,
thus suitable for a SOCP formulation.
So, in order to minimize the overall $\Delta v$, we can 
equivalently minimize the sum of the newly defined variables:
\begin{equation}
    J = \sum_{j = 1}^{M}{k^2 \rho_j \sigma_j}
    \label{eq:objective_linear}
\end{equation}

The resulting SOCP problem is:
\begin{align}
    \min_{\bm{\tilde{x}}_j,\, \bm{\Delta \tilde{v}}_j,\, \sigma_j}&
    J = \sum_{j = 1}^{M} k^2 \rho_j \sigma_j &
    \label{eq:objective_cvx} \\
    \text{subject to} 
        &\bm{\tilde{x}}_1^- = \bm{\tilde{x}}_0 & 
        \label{eq:initial_conditions_cvx}\\
        &\bm{\tilde{x}}_{M}^+ = \bm{\tilde{x}}_f & 
        \label{eq:final_conditions_cvx}\\
        &\bm{\tilde{x}}_{j+1}^- = \Phi(\theta_{j+1},\theta_j) \, 
        \bm{\tilde{x}}_j^+ 
        &&j = 1, \dots, M-1 
        \label{eq:differential_constraints_cvx} \\
        &\Delta \tilde{v}_{x, j}^2 + \Delta \tilde{v}_{y, j}^2 
        + \Delta \tilde{v}_{z, j}^2 \leq \sigma_j^2
        &&j = 1, \dots, M 
        \label{eq:cone_constraints_cvx} \\
        &\bm{\tilde{x}}_{j}^+ = \bm{\tilde{x}}_j^- 
        + B \, \bm{\Delta \tilde{v}}_j 
        &&j = 1, \dots, M
        \label{eq:before_after_burns_cvx}
\end{align}

\section{Numerical Results}
\label{sec:numerical_results}


In order to demonstrate the effectiveness of the proposed approach,
three test cases, whose optimal solution is available in literature, 
are here considered:
i) a four-impulse rendezvous with a target in a circular orbit,
ii) a three-impulse approach maneuver for the Automated Transfer Vehicle (ATV)
to the International Space Station, 
and iii) a reconfiguration of a spacecraft formation flying on a highly
elliptic orbit for the SIMBOL-X mission.
For the sake of simplicity, only the in-plane motion is shown here, 
as in-plane and out-plane motions are decoupled in the considered 
dynamical model.
Hereafter, the state variables $y$ and  $v_y$ are thus dropped.
Numerical results are provided for a uniform mesh with $M = 257$ nodes. 
The effect of the mesh size on the quality of the attained solutions 
is then discussed.

\subsection{Circle-to-Circle Rendezvous}
Let us consider a chaser spacecraft flying on a circular orbit of radius $r_0$
that has to rendezvous with a target on a circular orbit of radius 
$r_f = 1.2 r_0$. 
In particular, by setting the target orbit radius $r_f$ and mean motion equal to one, the problem boundary conditions are 
$\bm{r}_0 = \begin{bmatrix}  -\pi & 0 & 1/6 \end{bmatrix}$, 
$\bm{v}_0 = \begin{bmatrix} 1/4 & 0 & 0 \end{bmatrix}$, 
$\bm{r}_f = \begin{bmatrix} 0 & 0 & 0 \end{bmatrix}$, and 
$\bfv_f = \begin{bmatrix} 0 & 0 & 0 \end{bmatrix}$, respectively.
The mission must be accomplished in a limited time, corresponding to 
$\theta_f = $ \SI{10}{\radian}.

\begin{table}[hbt!]
    \renewcommand{\arraystretch}{1}
    \caption{Results for the four-impulse rendezvous}   
    \label{tab:four_impulse_results}
    \centering
    \begin{tabular}{c c c}
\hline
Parameter & Convex & Optimal \\ 
\hline
$\theta_1$ (rad) & 0 & 0 \\ 
$\theta_2$ (rad) & 2.8125 & 2.8033 \\ 
$\theta_3$ (rad) & 7.1875 & 7.1967 \\ 
$\theta_4$ (rad) & 10.0000 & 10.0000 \\ 
$\bm{\Delta v}_1$ & $\begin{bmatrix} -0.01596 & +0.00421 \end{bmatrix}$ & $\begin{bmatrix} -0.01575 & +0.00415  \end{bmatrix}$\\ 
$\bm{\Delta v}_2$ & $\begin{bmatrix} -0.02989 & +0.00162 \end{bmatrix}$ & $\begin{bmatrix} -0.03028 & +0.00158  \end{bmatrix}$\\ 
$\bm{\Delta v}_3$ & $\begin{bmatrix} +0.06309 & +0.00343 \end{bmatrix}$ & $\begin{bmatrix} +0.06387 & +0.00333  \end{bmatrix}$\\ 
$\bm{\Delta v}_4$ & $\begin{bmatrix} +0.06553 & +0.01719 \end{bmatrix}$ & $\begin{bmatrix} +0.06549 & +0.01724  \end{bmatrix}$\\ 
Overall $\Delta v$ & 0.17828 & 0.17828 \\ 
\hline
\end{tabular}

\end{table}

\begin{figure}[hbt!]
    \centering
    \subfigure[Inertial frame]{\label{fig:four_prussing_traj_eci}
        \includegraphics[width=0.45\linewidth]{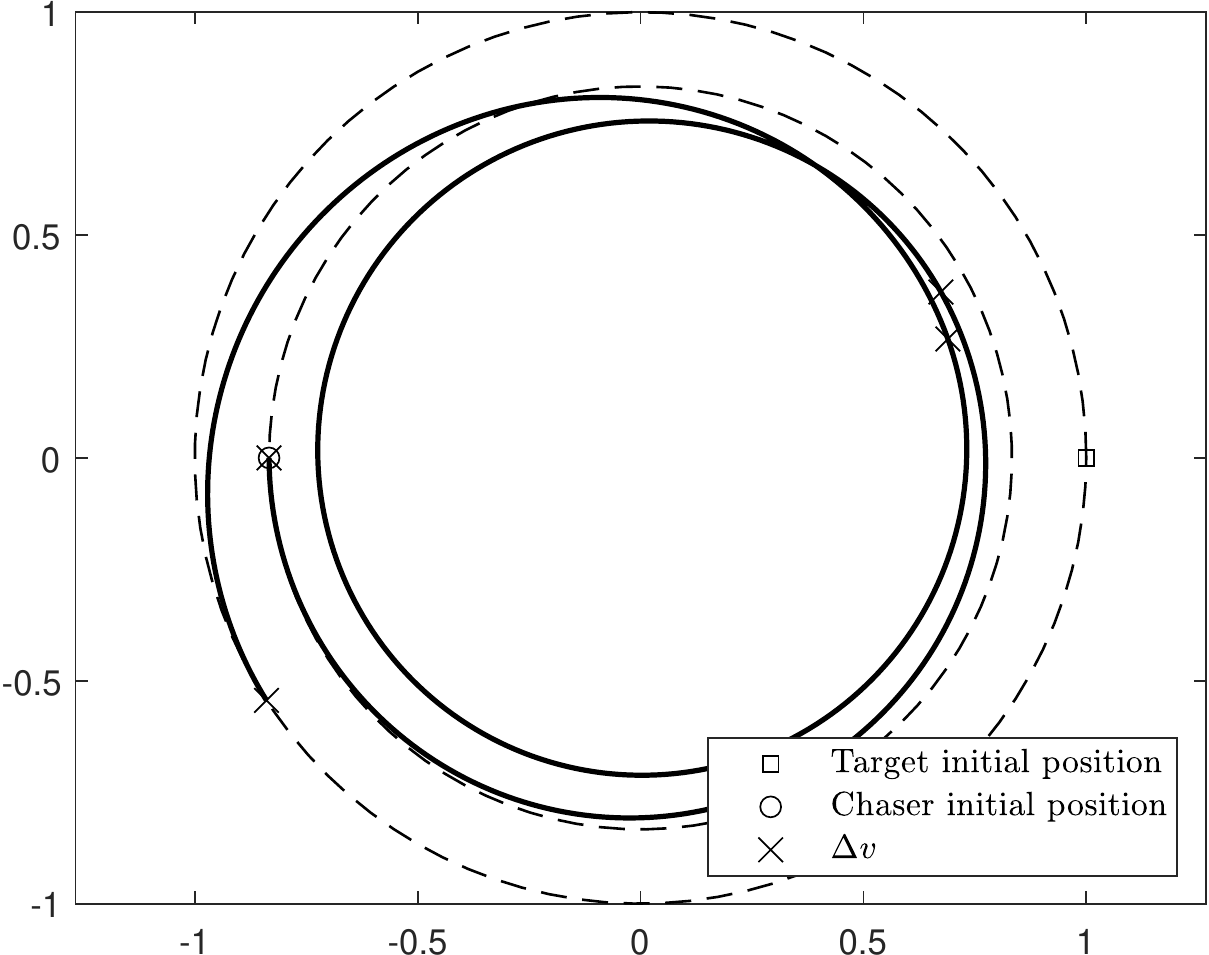}}
    \subfigure[Target orbital frame]{\label{fig:four_prussing_traj_hill}
        \includegraphics[width=0.45\linewidth]{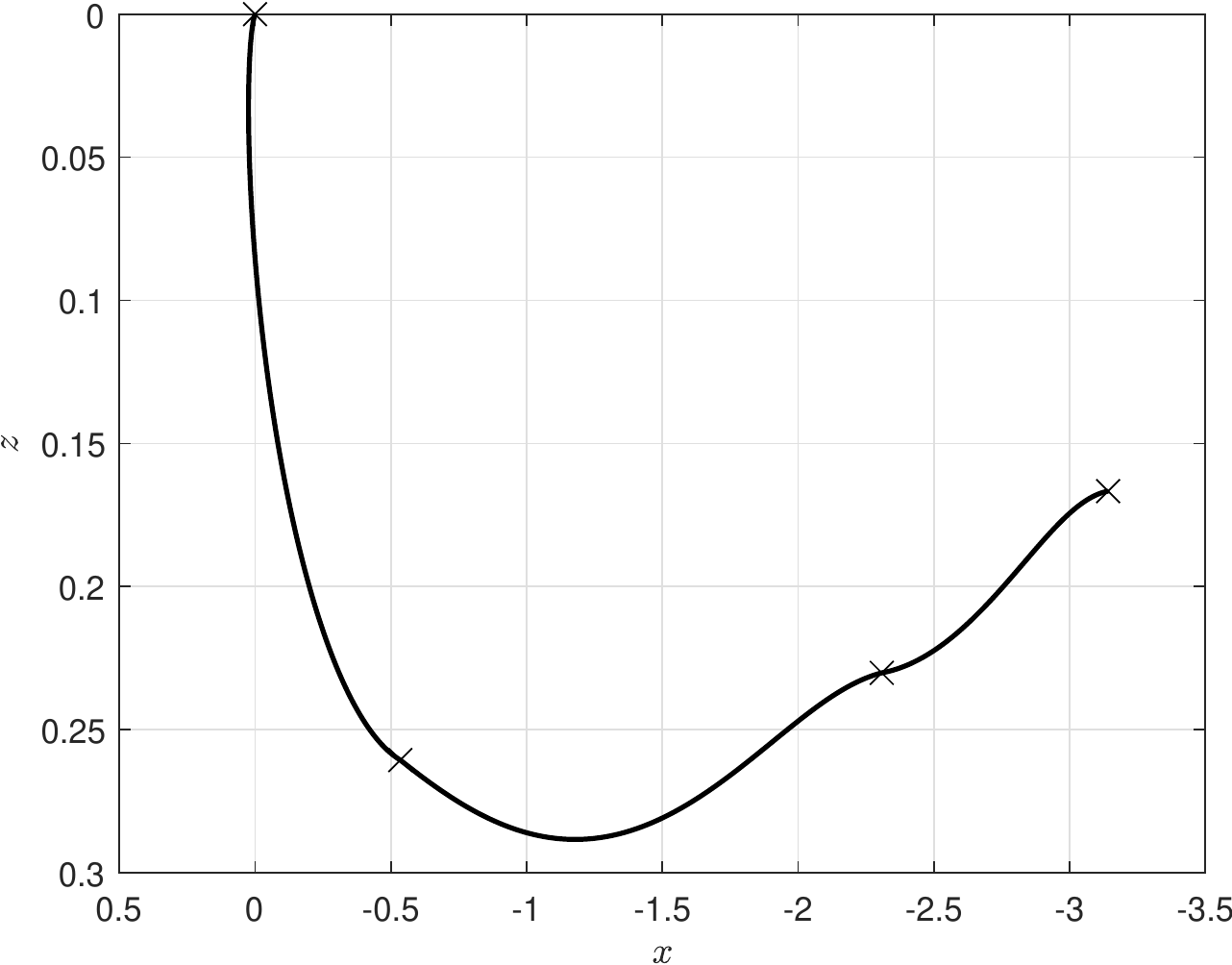}}
    \caption{Four-impulse trajectory}
   \label{fig:four_prussing_traj}
\end{figure}

Table~\ref{tab:four_impulse_results} compares the solution found by the 
present procedure with the optimal one, obtained by using an indirect 
technique that the authors routinely applied in the past for solving optimal 
control problems~\cite{Colasurdo1994, ZavoliAlaska, Simeoni2012}.
The convex approach succeeds in capturing the optimal 4-impulse structure of 
the mission, with one impulse at the departure, one at the arrival, and two 
burns at intermediate maneuvering points.
Besides these four points, the magnitude of the velocity changes $\sigma_j$
over the mesh points is several orders of magnitude lower than the 
actual impulses.
The two solutions are in good agreement. 
The overall $\Delta v$  is evaluated accurately. 
Minor differences arise in the locations of the intermediate impulses, 
which are captured at the grid points that are closer to the optimal locations 
found by the indirect method (as the adopted grid made up of 257 equally 
spaced nodes does not include the optimal $\theta_2$ and $\theta_3$). 
The shift of the intermediate impulses induces minimal changes in each 
$\bm{\Delta v}$ vector which, in turns, cause negligible changes in the 
overall trajectory (see Fig.~\ref{fig:four_prussing_traj}).

\subsection{ATV Approaching Maneuver}
The second case study concerns a more practical scenario, where the 
Automated Transfer Vehicle (ATV) needs to approach the International Space 
Station~\cite{AMADIEU199976}. 
An example is discussed by Louembet~\cite{louembet2017contributions}, 
where the optimal trajectory to move the ATV from about \SI{30}{\kilo\m} 
to \SI{100}{\m} to the ISS is investigated. 
More precisely, the following initial
conditions $\bfr_0 = \begin{bmatrix} -30 & 0 & 0.5   \end{bmatrix}$~km,  
$ \bfv_0 = \begin{bmatrix} 8.514 & 0 & 0 \end{bmatrix}$~m/s, 
and final conditions
$\bm{r}_f = \begin{bmatrix} -0.1 & 0 & 0    \end{bmatrix}$~km, 
$\bm{v}_f = \begin{bmatrix} 0 & 0 & 0 \end{bmatrix}$~m/s are prescribed.
%
The orbital elements of the target orbit and the transfer duration $\Delta t$ 
are reported in Table~\ref{tab:data_atv}.

\begin{table}[hbt!]
    \renewcommand{\arraystretch}{1}
    \caption{Data for the ATV rendezvous}
    \label{tab:data_atv}
    \centering
    \begin{tabular}{c c} 
\hline 
Parameter & Value \\ 
\hline 
$a$ (km) & \num{6763.00} \\ 
$e$ & \num{0.0052} \\ 
$i$ (deg) & \num{52.00} \\ 
$\Omega$ (deg) & \num{0.00} \\ 
$\omega$ (deg) & \num{0.00} \\ 
$\theta_0$ (deg) & \num{0.00} \\ 
$\Delta t$ (s) & \num{55350} \\ 
\hline 
\end{tabular} 

\end{table}

\begin{table}[hbt!]
    \renewcommand{\arraystretch}{1}
    \caption{Results for the ATV mission}
    \label{tab:results_atv}
    \centering
    \begin{tabular}{c c c c}
\hline
Parameter & Convex ($M = 257$) & Convex ($M^{\star} = 3$) & Optimal \\ 
\hline
$\theta_1$ (rad) & 0 & 0 & 0 \\ 
$\theta_2$ (rad) & 59.88620 & 59.90800 & 59.89691 \\ 
$\theta_3$ (rad) & 60.13170 &  &  \\ 
$\theta_4$ (rad) & 62.83150 & 62.83150 & 62.83149 \\ 
$\bm{\Delta v}_1$ (m/s) & $\begin{bmatrix} -7.55410 & +0.23649 \end{bmatrix}$ & $\begin{bmatrix} -7.55410 & +0.23952 \end{bmatrix}$ & $\begin{bmatrix} -7.55410 & +0.23877 \end{bmatrix}$\\ 
$\bm{\Delta v}_2$ (m/s) & $\begin{bmatrix} +0.13839 & +0.00141 \end{bmatrix}$ & $\begin{bmatrix} +0.14439 & +0.00085 \end{bmatrix}$ & $\begin{bmatrix} +0.14422 & +0.00086 \end{bmatrix}$\\ 
$\bm{\Delta v}_3$ (m/s) & $\begin{bmatrix} +0.00564 & +0.00003 \end{bmatrix}$ &  & \\ 
$\bm{\Delta v}_4$ (m/s) & $\begin{bmatrix} +0.04138 & +0.00126 \end{bmatrix}$ & $\begin{bmatrix} +0.04125 & +0.00131 \end{bmatrix}$ & $\begin{bmatrix} +0.04147 & +0.00130 \end{bmatrix}$\\ 
Overall $\Delta v$ (m/s) & 7.74357 & 7.74356 & 7.74356 \\ 
\hline
\end{tabular}

\end{table}

The optimal solution and the one attained for a $257$-node grid are reported 
in Table~\ref{tab:results_atv}.
The optimal transfer found out by 
Louembet~\cite{louembet2017contributions},
is made up of three impulses: one at the departure, one at the arrival, 
and one at an intermediate point.
This problem is numerically more challenging than the previous one, 
as the optimal solution presents two burns much smaller than the initial one, 
and the change in the transversal direction is of about 
three orders of magnitude.
Indeed, the convex approach shows some difficulties in determining
precisely the exact location of the intermediate maneuver. 
Rather, the intermediate impulse is spread over adjacent nodes.
Eventually, a "sub-optimal" four-impulse solution is detected:
the \nth{2} impulse is close to the optimal one, 
the \nth{3} impulse is much lower but it is retained, 
because greater than the allowed tolerance $10^{-5}$.
However, the difference in the overall mission cost is about 
$10^{-5}$, which is close to the optimization tolerance.

For the sake of completeness, a parametric investigation has been carried out, 
by considering a 3-nodes grid and moving the inner node location over the 
domain.
A clear minimum is identified and the best-found solution (named $M^{\star}=3$)
is also reported in Table~\ref{tab:results_atv}.

\subsection{SIMBOL-X}
Finally, let us consider a case where the reference orbit is highly elliptical,
as in the case of the SIMBOL-X mission~\cite{gaudel2010autonomous}.
An approach maneuver from \SI{30}{\kilo\m} to \SI{500}{\m} 
from the target is investigated.
Specifically, the initial
conditions $\bfr_0 = \begin{bmatrix} 18.3095 & 0 & -23.7647 \end{bmatrix}$ km,
$\bfv_0 = \begin{bmatrix} -0.0542 & 0 & -0.0418 \end{bmatrix}$ m/s, 
and final conditions
$\bm{r}_f = \begin{bmatrix} 335.12 & 0 & -371.1 \end{bmatrix} \text{m}$, 
$\bm{v}_f = \begin{bmatrix} 0.00155 & 0 & 0.0014 \end{bmatrix} \text{m/s}$  
are prescribed.
%
%
Table~\ref{tab:data_simbol_x} reports the target orbital elements 
and the duration of the maneuver.

\begin{table}[hbt!]
    \renewcommand{\arraystretch}{1}
    \caption{Data for the SIMBOL-X rendezvous}
    \label{tab:data_simbol_x}
    \centering
    \begin{tabular}{c c} 
\hline 
Parameter & Value \\ 
\hline 
$a$ (km) & \num{106246.98} \\ 
$e$ & \num{0.7988} \\ 
$i$ (deg) & \num{5.20} \\ 
$\Omega$ (deg) & \num{90.00} \\ 
$\omega$ (deg) & \num{180.00} \\ 
$\theta_0$ (deg) & \num{135.00} \\ 
$\Delta t$ (s) & \num{49995} \\ 
\hline 
\end{tabular} 

\end{table}

\begin{table}[hbt!]
    \renewcommand{\arraystretch}{1}
    \caption{Results for the SIMBOL-X mission}
    \label{tab:results_simbol_x}
    \centering
    \begin{tabular}{c c c c}
\hline
Parameter & Convex & Optimal & \\ \hline
$\theta_1$ & 2.3562 & 2.3562 \\ 
$\theta_2$ & 2.7859 & 2.7859 \\ 
$\bm{\Delta v}_1$ (m/s) & $\begin{bmatrix} -0.6193 & +0.5061 \end{bmatrix}$ & $\begin{bmatrix} -0.6193 & +0.5061 \end{bmatrix}$\\ 
$\bm{\Delta v}_2$ (m/s) & $\begin{bmatrix} +0.1748 & -0.4912 \end{bmatrix}$ & $\begin{bmatrix} +0.1748 & -0.4912 \end{bmatrix}$\\ 
Overall $\Delta v$ (m/s) & 1.3212 & 1.3212 \\ 
\hline
\end{tabular}

\end{table}

The final solution is a two-impulse transfer, as reported in 
Table~\ref{tab:results_simbol_x}.
The solution attained with the convex approach is identical 
to the optimal one, provided by Arzelier et al.~\cite{arzelier2011using}. 
Indeed, in this case, the impulses are located at the initial and final time, 
respectively, which belong to the mesh grid; 
hence, the impulses are captured at their exact optimal location.

\subsection{Effects of the Mesh Grid}
\begin{figure}[hbt!]
    \centering
    \includegraphics[width=.5\textwidth]{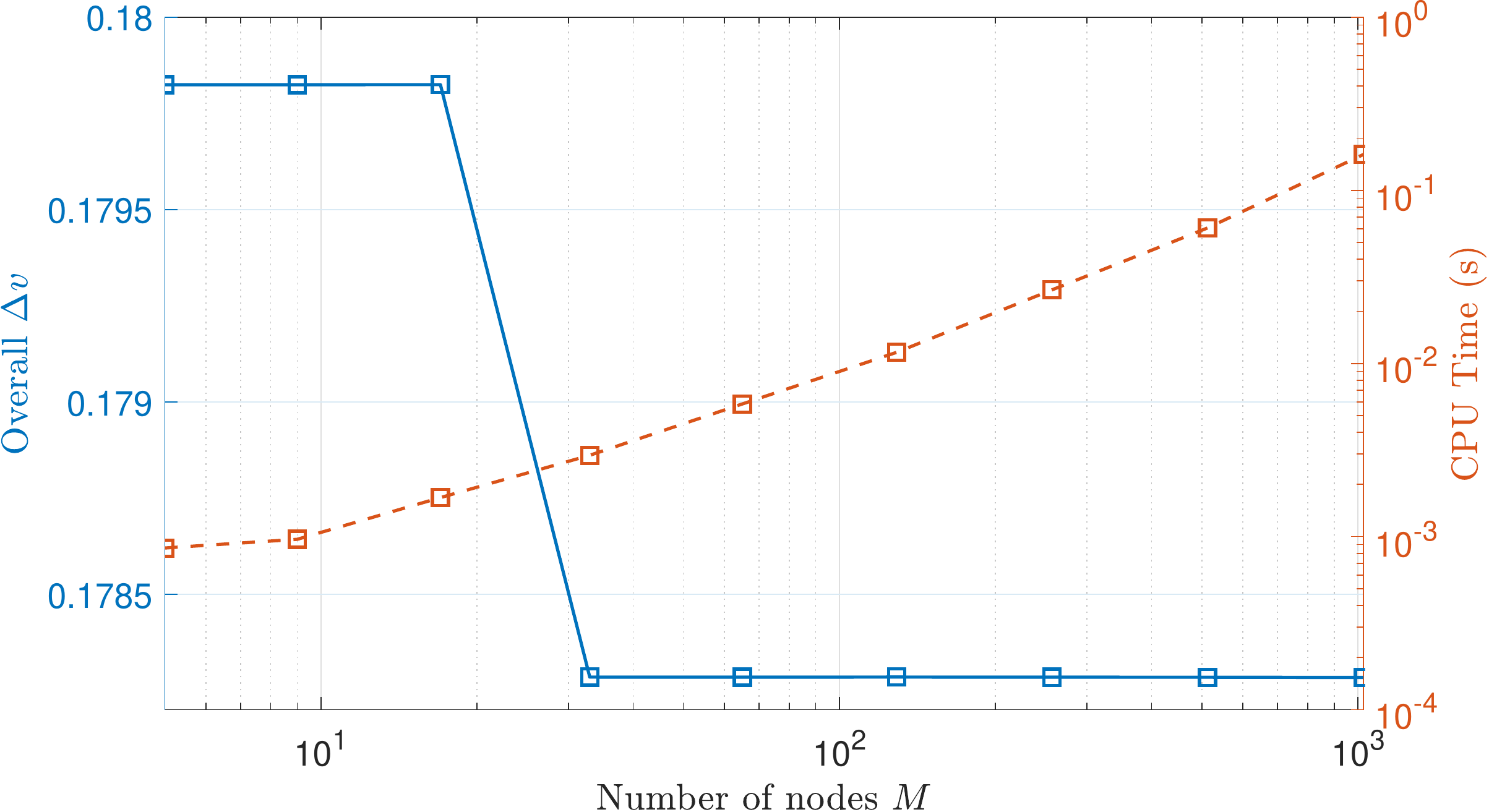}
    \caption{Four-impulse overall $\Delta v$ and computation time as a 
        function of the number of discretization nodes 
        for the circle-to-circle example}
    \label{fig:DV_and_CPU_time_prussing}
\end{figure}

A parametric study on the effects of the mesh size has also been carried out, 
in order to investigate its role on the solution quality in terms of 
overall $\Delta v$ estimation, computational time, and accuracy on the 
location of the maneuvering points.
Intuitively, one expects that, as the mesh density increases, 
the distance between the optimal locations of the internal impulses and 
the (closest) grid nodes decreases, 
leading to a numerical solution closer to the optimal one.
Figure~\ref{fig:DV_and_CPU_time_prussing} shows the overall mission cost 
for the circle-to-circle rendezvous mission,
and the corresponding computational time,
as a function of the number of grid points.
The overall $\Delta v$ quickly converges towards the optimal value.
Instead, the computational cost increases almost linearly in the number of 
nodes.
The use of a mesh grid with a small, limited, number of nodes seems thus 
justified.
It is also worth mentioning that 
the  computation time is quite short even for the largest attempted grid
because the convex formulation allows for the use of highly-efficient algorithms.

\begin{figure}[hbt!]
    \centering
    \includegraphics[width=\textwidth]{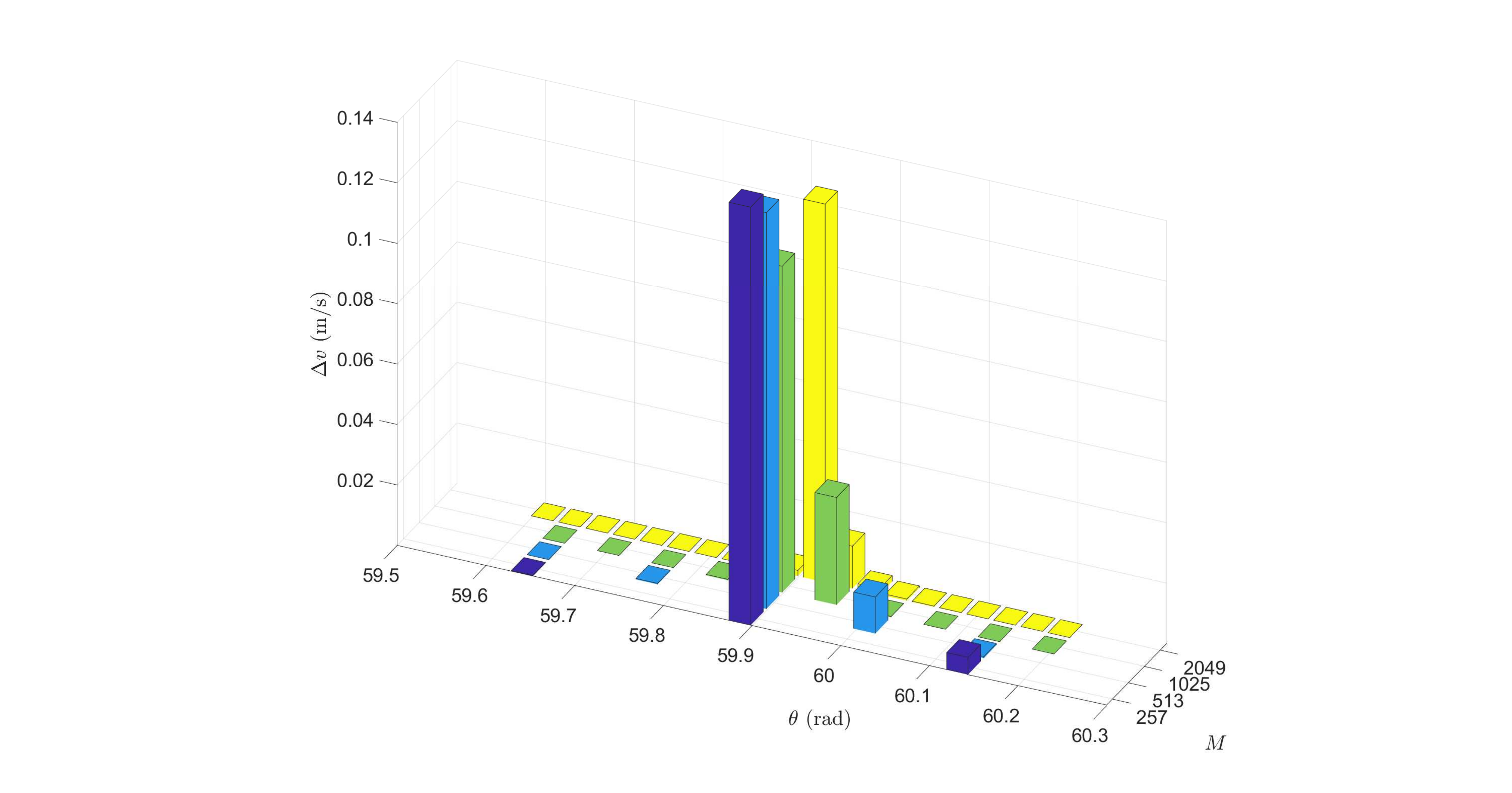}
    \caption{Distribution of the $\Delta v$ magnitudes for different mesh sizes
    in the neighborhood of the intermediate impulse in the ATV test case.}
    \label{fig:histogram_atv3}
\end{figure}

On the other hand,
increasing the number of grid points may lead to numerical issues related to 
the accumulation of truncation and round-off errors. 
Due to the finite precision (i.e., tolerance) of the optimization process, in 
some circumstances, an impulse may be spread over a set of 
neighboring 
nodes. 
An example is proposed in Figure~\ref{fig:histogram_atv3},
where the distribution of the velocity change magnitudes in the neighborhood 
of the (optimal) intermediate impulse location for the ATV test case is 
presented.
Even though this issue could be mitigated by tightening the optimization 
tolerances and by selecting a suitable set of nondimensionalization factors 
to properly scale the problem, 
the current analysis does not provide a systematic way to avoid this problem,
which manifests when too many grid nodes are located in the proximity of the 
optimal maneuvering point.

\section{Conclusion}
\label{sec:conclusion}

This note presented a convex formulation of the time-fixed optimal rendezvous in a linear dynamics. The Tschauner-Hempel equations are used to describe the motion of the chaser spacecraft relative to a target that flies on a keplerian orbit of arbitrary eccentricity.
The original non-convex problem is transformed into a convex one by using a 
gridding technique, 
that is, 
by introducing a finite discretization of the independent variable domain: 
impulse locations are constrained to the grid nodes, 
and the differential constraints are replaced by algebraic constraints 
involving the state transition matrix between adjacent nodes.
The convex formulation allows the use of special solution algorithms that 
guarantee the convergence towards the global optimum within a limited 
computation time, 
even when a large number of variables are involved, 
and that does not require any sort of initialization.

Numerical results show a good agreement between the convex solution and 
the optimal one, provided by other approaches.
In cases of practical interest,  
a limited number of nodes is sufficient to get an accurate solution in terms 
of both trajectory and mission cost. 
A parametric study on the effects of the mesh size discourages the use of 
very-fine grids.
In fact, the accumulation of truncation errors and the finite precision of the 
adopted convex solver may hinder the effectiveness of the approach.
Instead, the use of medium-to-coarse grids allows
to get the most out of this approach:
the overall $\Delta v$ cost quickly converges towards the optimal value, and 
the computational effort is kept to the minimum.
These features make the proposed approach suitable for time-critical 
applications, such as autonomous guidance, 
and other computationally demanding tasks.

\bibliography{bibliography}

\begin{thebibliography}{33}
\newcommand{\enquote}[1]{``#1''}
\providecommand{\natexlab}[1]{#1}
\providecommand{\url}[1]{\texttt{#1}}
\providecommand{\urlprefix}{URL }
\expandafter\ifx\csname urlstyle\endcsname\relax
  \providecommand{\doi}[1]{\discretionary{}{}{}https://doi.org/#1}\else
  \providecommand{\doi}[1]{\discretionary{}{}{}\urlstyle{rm}\url{https://doi.org/#1}}\fi

\bibitem[{Weiss et~al.(2015)Weiss, Baldwin, Erwin, and
  Kolmanovsky}]{weiss2015model}
Weiss, A., Baldwin, M., Erwin, R.~S., and Kolmanovsky, I., \enquote{Model
  predictive control for spacecraft rendezvous and docking: Strategies for
  handling constraints and case studies,} \emph{IEEE Transactions on Control
  Systems Technology}, Vol.~23, No.~4, 2015, pp. 1638--1647.

\bibitem[{Subbarao and Welsh(2008)}]{subbarao2008nonlinear}
Subbarao, K., and Welsh, S.~J., \enquote{Nonlinear control of motion
  synchronization for satellite proximity operations,} \emph{Journal of
  guidance, control, and dynamics}, Vol.~31, No.~5, 2008, pp. 1284--1294.

\bibitem[{Vassar and Sherwood(1985)}]{vassar1985formationkeeping}
Vassar, R.~H., and Sherwood, R.~B., \enquote{Formationkeeping for a Pair of
  Satellites in a Circular Orbit,} \emph{Journal of Guidance, Control, and
  Dynamics}, Vol.~8, No.~2, 1985, pp. 235--242.

\bibitem[{Epenoy(2011)}]{epenoy2011fuel}
Epenoy, R., \enquote{Fuel-optimal trajectories for continuous-thrust orbital
  rendezvous with collision avoidance constraint,} \emph{Advances in the
  Astronautical Sciences}, Vol. 140, No.~2, 2011, pp. 341--360.

\bibitem[{Zavoli et~al.(2019)Zavoli, Federici, Benedikter, Casalino, and
  Colasurdo}]{zavoli2019gtoc}
Zavoli, A., Federici, L., Benedikter, B., Casalino, L., and Colasurdo, G.,
  \enquote{GTOC X: Solution Approach of Team Sapienza-PoliTo,} \emph{Paper AAS
  19-894, Astrodynamics Specialist Conference}, Portland, Maine, 2019.

\bibitem[{Hill(1878)}]{hill1878researches}
Hill, G.~W., \enquote{Researches in the lunar theory,} \emph{American journal
  of Mathematics}, Vol.~1, No.~1, 1878, pp. 5--26.

\bibitem[{Baranov(1990)}]{Baranov1990}
Baranov, A., \enquote{An algorithm for calculating parameters of multi-orbit
  maneuvers in remote guidance,} \emph{Cosmic Research}, Vol.~28, No.~1, 1990,
  pp. 61--67.

\bibitem[{Tschauner and Hempel(1965)}]{tschauner1965rendezvous}
Tschauner, J., and Hempel, P., \enquote{Rendezvous zu einem in elliptischer
  Bahn umlaufenden Ziel,} \emph{Astronautica Acta}, Vol.~11, No.~2, 1965, pp.
  104--+.

\bibitem[{Yamanaka and Ankersen(2002)}]{Yamanaka2002}
Yamanaka, K., and Ankersen, F., \enquote{New State Transition Matrix for
  Relative Motion on an Arbitrary Elliptical Orbit,} \emph{Journal of Guidance,
  Control, and Dynamics}, Vol.~25, No.~1, 2002, pp. 60--66.
\newblock \doi{10.2514/2.4875}.

\bibitem[{Kechichian(1992)}]{Kechichian1992}
Kechichian, J., \enquote{Techniques of accurate analytic terminal rendezvous in
  near-circular orbit,} \emph{Acta Astronautica}, Vol.~26, No.~6, 1992, pp.
  377--394.
\newblock \doi{10.1016/0094-5765(92)90068-t}.

\bibitem[{Schweighart and Sedwick(2002)}]{Schweighart2002}
Schweighart, S.~A., and Sedwick, R.~J., \enquote{High-Fidelity Linearized J
  Model for Satellite Formation Flight,} \emph{Journal of Guidance, Control,
  and Dynamics}, Vol.~25, No.~6, 2002, pp. 1073--1080.
\newblock \doi{10.2514/2.4986}.

\bibitem[{Mirfakhraie and Conway(1994)}]{Mirfakhraie1994}
Mirfakhraie, K., and Conway, B.~A., \enquote{Optimal cooperative time-fixed
  impulsive rendezvous,} \emph{Journal of Guidance Control Dynamics}, Vol.~17,
  1994, pp. 607--613.
\newblock \doi{10.2514/3.21240}.

\bibitem[{Lu and Liu(2013)}]{lu2013autonomous}
Lu, P., and Liu, X., \enquote{Autonomous trajectory planning for rendezvous and
  proximity operations by conic optimization,} \emph{Journal of Guidance,
  Control, and Dynamics}, Vol.~36, No.~2, 2013, pp. 375--389.

\bibitem[{Benedikter et~al.(2019)Benedikter, Zavoli, and
  Colasurdo}]{benedikter2019convexrendezvous}
Benedikter, B., Zavoli, A., and Colasurdo, G., \enquote{A Convex Optimization
  Approach for Finite-Thrust Time-Constrained Cooperative Rendezvous,}
  \emph{Paper AAS 19-763, Astrodynamics Specialist Conference}, Portland,
  Maine, 2019.

\bibitem[{Zavoli and Colasurdo(2015)}]{Zavoli2015}
Zavoli, A., and Colasurdo, G., \enquote{Indirect optimization of finite-thrust
  cooperative rendezvous,} \emph{Journal of Guidance, Control, and Dynamics},
  Vol.~38, No.~2, 2015, pp. 304--314.
\newblock \doi{10.2514/1.G000531}, cited By 6.

\bibitem[{Neustadt(1964)}]{neustadt1964optimization}
Neustadt, L.~W., \enquote{Optimization, a moment problem, and nonlinear
  programming,} \emph{Journal of the Society for Industrial and Applied
  Mathematics, Series A: Control}, Vol.~2, No.~1, 1964, pp. 33--53.

\bibitem[{Lawden(1963)}]{lawden1963general}
Lawden, D., \enquote{General theory of optimal rocket trajectories,}
  \emph{Optimal trajectories for space navigation}, 1963, pp. 54--78.

\bibitem[{Prussing(1994)}]{prussing1994optimal}
Prussing, J.~E., \enquote{Optimal impulsive linear systems- Sufficient
  conditions and maximum number of impulses,} \emph{Spaceflight mechanics
  1994}, 1994, pp. 1015--1023.

\bibitem[{Carter(1991)}]{carter1991optimal}
Carter, T., \enquote{Optimal impulsive space trajectories based on linear
  equations,} \emph{Journal of optimization theory and applications}, Vol.~70,
  No.~2, 1991, pp. 277--297.

\bibitem[{Jezewski(1980)}]{Jezewski1980}
Jezewski, D., \enquote{Primer vector theory applied to the linear
  relative-motion equations,} \emph{Optimal Control Applications and Methods},
  Vol.~1, No.~4, 1980, pp. 387--401.
\newblock \doi{10.1002/oca.4660010408}.

\bibitem[{Prussing(1969)}]{prussing1969optimal}
Prussing, J., \enquote{Optimal four-impulse fixed-time rendezvous in the
  vicinity of a circular orbit.} \emph{AIAA Journal}, Vol.~7, No.~5, 1969, pp.
  928--935.

\bibitem[{Carter and Alvarez(2000)}]{carter2000quadratic}
Carter, T.~E., and Alvarez, S.~A., \enquote{Quadratic-based computation of
  four-impulse optimal rendezvous near circular orbit,} \emph{Journal of
  Guidance, Control, and Dynamics}, Vol.~23, No.~1, 2000, pp. 109--117.

\bibitem[{LION and HANDELSMAN(1968)}]{Lion1968}
LION, P.~M., and HANDELSMAN, M., \enquote{Primer vector on fixed-time impulsive
  trajectories.} \emph{{AIAA} Journal}, Vol.~6, No.~1, 1968, pp. 127--132.
\newblock \doi{10.2514/3.4452}.

\bibitem[{Chernick and D’Amico(2017)}]{chernick2017new}
Chernick, M., and D’Amico, S., \enquote{New closed-form solutions for optimal
  impulsive control of spacecraft relative motion,} \emph{Journal of Guidance,
  Control, and Dynamics}, Vol.~41, No.~2, 2017, pp. 301--319.

\bibitem[{Prussing(2003)}]{prussing2003optimal}
Prussing, J.~E., \enquote{Optimal two-and three-impulse fixed-time rendezvous
  in the vicinity of a circular orbit,} \emph{Journal of Spacecraft and
  Rockets}, Vol.~40, No.~6, 2003, pp. 952--959.

\bibitem[{Alizadeh and Goldfarb(2003)}]{alizadeh2003second}
Alizadeh, F., and Goldfarb, D., \enquote{Second-order cone programming,}
  \emph{Mathematical programming}, Vol.~95, No.~1, 2003, pp. 3--51.

\bibitem[{Colasurdo and Pastrone(1994)}]{Colasurdo1994}
Colasurdo, G., and Pastrone, D., \enquote{Indirect optimization method for
  impulsive transfers,} \emph{Astrodynamics Conference}, American Institute of
  Aeronautics and Astronautics, 1994.
\newblock \doi{10.2514/6.1994-3762}.

\bibitem[{Zavoli et~al.(2012)Zavoli, Simeoni, Casalino, and
  Colasurdo}]{ZavoliAlaska}
Zavoli, A., Simeoni, F., Casalino, L., and Colasurdo, G., \enquote{Optimal
  cooperative deployment of a two-satellite formation into a highly elliptic
  orbit,} \emph{Advances in the Astronautical Sciences}, Vol. 142, 2012, pp.
  3647--3663.

\bibitem[{Simeoni et~al.(2012)Simeoni, Casalino, Zavoli, and
  Colasurdo}]{Simeoni2012}
Simeoni, F., Casalino, L., Zavoli, A., and Colasurdo, G., \enquote{Indirect
  optimization of satellite deployment into a highly elliptic orbit,}
  \emph{International Journal of Aerospace Engineering}, Vol. 2012, 2012,
  p.~14.
\newblock \doi{10.1155/2012/152683}.

\bibitem[{Amadieu and Heloret(1999)}]{AMADIEU199976}
Amadieu, P., and Heloret, J., \enquote{The automated transfer vehicle,}
  \emph{Air and Space Europe}, Vol.~1, No.~1, 1999, pp. 76 -- 80.
\newblock \doi{https://doi.org/10.1016/S1290-0958(99)80044-6},
  \urlprefix\url{http://www.sciencedirect.com/science/article/pii/S1290095899800446}.

\bibitem[{Louembet(2017)}]{louembet2017contributions}
Louembet, C., \enquote{Contributions au guidage pour le rendez-vous spatial par
  r{\'e}solution du probl{\`e}me de commande optimale impulsionnelle,} Ph.D.
  thesis, Universit{\'e} Toulouse 3 Paul Sabatier (UT3 Paul Sabatier), 2017.

\bibitem[{Gaudel et~al.(2010)Gaudel, Berges, Trapier, Gamet, and
  Djalal}]{gaudel2010autonomous}
Gaudel, A., Berges, J.-C., Trapier, T., Gamet, P., and Djalal, S.,
  \enquote{Autonomous rendezvous guidance function of the SIMBOL-X formation
  flying mission a high elliptical orbit: Preliminary design and performance
  analysis,} \emph{Proceedings of the 21st International Symposium on Space
  Flight Dynamics}, 2010.

\bibitem[{Arzelier et~al.(2011)Arzelier, Kara-Zaitri, Louembet, and
  Delibasi}]{arzelier2011using}
Arzelier, D., Kara-Zaitri, M., Louembet, C., and Delibasi, A., \enquote{Using
  polynomial optimization to solve the fuel-optimal linear impulsive rendezvous
  problem,} \emph{Journal of Guidance, Control, and Dynamics}, Vol.~34, No.~5,
  2011, pp. 1567--1576.

\end{thebibliography}

\end{document}